\pdfoutput=1
\documentclass[%
aps,								%sigla della rivista: APS, nello specifico Physical Review Letters
prl,								%N.B.: nella modalità prl NON VENGONO NUMERATE lE SEZIONI!!!
reprint,							  %opzioni per 1 o 2 colonne, variabile a seconda della rivista
showpacs, showkeys,		  	%opzioni per i codici e le parole chiave dell'articolo
groupedaddress,					  %opzioni per gli indirizzi degli autori
nobibnotes,		    %SCOMMENTATO
floatfix,				 %SCOMMENTATO
longbibliography,%AGGIUNTA MIA
]{revtex4-1}
\usepackage{dcolumn}							% Align table columns on decimal point
\usepackage[mathlines]{lineno}% Enable numbering of text and display math
\usepackage{color}
\usepackage{placeins}
\usepackage{tabularx,colortbl}
\usepackage{ifpdf}
\ifpdf			%COMMENTATO
   \usepackage[pdftex]{graphicx}
\else								%COMMENTATO
   \usepackage[dvips]{graphicx}
\fi									%COMMENTATO
\usepackage{hyperref}
\hypersetup{
pdfstartview={XYZ},
pdftitle = {Cross_product_ND_prl},
pdfauthor = {Carlo Andrea Gonano}
}

\usepackage{comment}
\usepackage{url}									%Pacchetto per riferimenti ai siti internet
\usepackage{natbib}								%Pacchetto utile per la bibliografia di BibTeX
\usepackage{textcase}
\usepackage{hypernat}	
\usepackage{dcolumn} 							%pacchetto per allineare i punti alla virgola
\usepackage{longtable} 						%pacchetto per la gestione delle tabelle lunghe
\usepackage{amsfonts}
\usepackage{amsmath} %fornisce svariate estensioni per il miglioramento della gestione
\usepackage{amssymb} %arricchisce la scelta di simboli matematici
\usepackage{amsthm}  %migliora la gestione degli enunciati matematici,
\usepackage{mathtools}
\usepackage{bm} %Pacchetto per il corsivo neretto matematico
\usepackage{empheq} %pacchetto per evidenziare formule ed equazioni.
\newcommand{\virga}{``}
\newcommand{\virgc}{''\,}
\DeclarePairedDelimiter{\norma}{\lVert}{\rVert}

\newcommand{\vect}{\vec}	%normale freccia
\newcommand{\vectnabla}{\overrightarrow{\nabla}} 	%normale freccia
\newcommand{\prodscal}[2]%
{\vect{#1}^T \!\!\cdot\!\vect{#2} \;}

\newcommand{\prodvet}[2]%
{\vect{#1} \!\wedge \!\vect{#2} \;}

\newcommand{\crossprod}[2]%
{\vect{#1} \!\times \!\vect{#2} \;}

\newcommand{\triprodscal}[3]%
{ (\prodvet{#1}{#2})^T \!\! \cdot \! \vect{#3} }

\newcommand{\doubleoverline}[1]%
{\overline{\overline{{#1}}}}

\newcommand{\proddiad}[2]%
{[\overline{\overline{{#1}{#2 ^T}}} ]}
\newcommand{\doublewedge}{ \: \overset{\wedge} {\text{\scriptsize{$\wedge$}}} \; }%

\newcommand{\prodvett}[2]%
{\vect{#1} \doublewedge \vect{#2} }

\newcommand{\prodvettdiad}[2]%
{\proddiad{#2}{#1}  - \proddiad{#1}{#2}}
\newcommand{\gradmat}[2]%
{
\frac{\partial {#1}}{\partial {#2} }
}
\newcommand{\gradmatt}[2]%
{\left[
\frac{\partial {#1}}{\partial {#2} }
\right]}
\newcommand{\divergv}[1]
{\vectnabla^T \!\!\cdot \!\vect{#1}}

\newcommand{\rotor}[1]
{\vectnabla \wedge \vect{#1}}
\newcommand{\rotorr}[1]
{\vectnabla \doublewedge \vect{#1}}
\newcommand{\rotorrdiad}[2]
{\gradmatt{\vect{#1}}{\vect{#2}} - \gradmatt{\vect{#1}}{\vect{#2}}^T}

\begin{document}
%
%===============================================================
%
%											*** FRONTMATTER ***
%
%===============================================================
%
%\preprint{APS/123-QED}
\title{Cross product in N Dimensions - the doublewedge product}
%use linebreaks \\ within to get better formatting as desired
%-------------------------------------------------------------------
% VERSIONE DI REVTEX4.1
%\preprint{}

%AUTORI
%
\author{Carlo Andrea \surname{Gonano}}%
\email[e-mail: ]{carloandrea.gonano@polimi.it}
%\thanks{carloandrea.gonano@polimi.it bis}
\author{Riccardo Enrico \surname{Zich}}%
\email[e-mail: ]{riccardo.zich@polimi.it}
%\thanks{}
\affiliation{Politecnico di Milano, Energy Department, 
via La Masa 34, 20156 Milan, MI, Italy} %<-- affiiliazione comune
% In caso di molti autori, usare superscriptaddress in \documentclass
%
% N.B.: nella modalità prl (Physical Review Letters) della \documentclass
% gli indirizzi e-mail sono spostati appena prima della bibliografia,
% anziché al fondo della 1a pagina. Per riportarli alla posizione giusta,
% usare l'opzione groupedaddress in \documentclass{}
%
%------------------------------------------------------------------
%DATA
\date{\today}
\begin{abstract}
%
%NORMA ICEAA
%The paper should start with an abstract giving a short overview on the discussed matter and the presented results. The %abstract should be set using 8pt, bold face font. It should not exceed 15 lines.
%
%NORMA IEEE
%\boldmath
%Please include a brief abstract here. The abstract should be limited to 50-200 words and should concisely state what was 
%done, how it was done, principal results, and their significance.

%NORMA PHYSICAL REVIEW
%An abstract must accompany each manuscript. The abstract should consist of one paragraph and be com-
%pletely self-contained. It cannot contain numbered references; incorporate such information into the
%abstract itself. Use this form:
%Further information is available [A. B. Smith, Phys. Rev. A
%26, 107 (1982)].
%Displayed equations and tabular material are discour aged.
%Define all nonstandard symbols and abbreviations.
%------------------------------------------------------------------------
%VERSIONE INTEGRALE
% 
%\begin{comment}
The cross product $\times$ frequently occurs in Physics and Engineering, since it has large applications in many contexts, 
e.g. for calculating angular momenta, torques, rotations, volumes etc.
Though this mathematical operator is widely used, it is commonly expressed in a 3-D notation which gives rise to many 
paradoxes and difficulties. In fact, instead of other vector operators like scalar product, the cross product is defined just 
in 3-D space, it does not respect reflection rules and invokes the concept of \virga handedness\virgc.
%In this paper we are going to present an extension of cross product in N spatial Dimensions,
In this paper we are going to present an extension of cross product in an arbitrary number N of spatial Dimensions,
different from the one adopted in 
the Exterior Algebra and explicitly designed for an easy calculus of moments.  
%\end{comment}
%
%120 PAROLE, 8 righe centrate
%6** caratteri (senza spazi)
%7** caratteri (con spazi)
%
%
% VERSIONE SUPER-RIDOTTA
\begin{comment}
The cross product frequently occurs in Physics and Engineering, since it has large applications in many contexts, e.g. for
calculating angular momenta, torques, rotations, volumes etc. Though this mathematical operator is widely used, it is
commonly expressed in a 3-D notation which gives rise to many paradoxes and difficulties. In fact, instead of scalar
product, the cross product is defined just in 3-D space and it does not respect reflection rules. In this paper we present
an extension of cross product in N spatial Dimensions, explicitly designed for an easy calculus of moments.
\end{comment}
%    
%
%93 PAROLE, 6 righe centrate
%496 caratteri (senza spazi)
%588 caratteri (con spazi)
%
%
%
%
%----------------------------------------------------------------------------------
\end{abstract}
%
%-------------------
% Norme IEEE
% AGGIUNTE DI REVTEX4.1
%---------------------------------------------
%
%PACS = PHYSICAL AND ASTRONOMY CLASSIFICATION SCHEME
%
%Each manuscript must be assigned indexing codes
%which are used in computerized secondary information
%services. See also Physical Review Letters, 14 December
%1992, for code indexing information. In general, follow
%these guidelines.
%(1) Choose no more than four index number codes.
%(2) Place your principal index code  first.
%(3) Always choose the lowest-level code available.
%(4) Always include the check characters.
%
%All indexing will be veried by the journal scientic editor.
%
%http://journals.aps.org/PACS/
%http://www.aip.org/pacs
%\pacs{23.23.+x, 56.65.Dy DA SCRIVERE}		%ESEMPIO DI CODICE PACS PER LA SOTTOMISSIONE AD APS
\pacs{45.20.-d, 45.10.Na, 02.40.Yy, 45.20.da}
%PACS= the Physics and Astronomy  Classification Scheme.
%----------------------------------------------------------------------------
%Per articolo su prodotto vettore N-D:
%
%00-General
%http://www.aip.org/publishing/pacs/pacs-reg00#02
%
% 01. Communication, education, history, and philosophy
% 02. Mathematical methods in physics
%		02.10.Ud	Linear algebra
%		02.10.Yn	Matrix theory
%		02.30.Tb	Operator theory
%		02.30.Vv	Operational calculus
%		02.40.Dr	Euclidean and projective geometries
%		02.40.Yy	Geometric mechanics (see also 45.20.Jj in formalisms in classical mechanics)
%
%40—Electromagnetism, Optics, Acoustics, Heat Transfer,
%Classical Mechanics, and Fluid Dynamics
% 45. Classical mechanics of discrete systems
%		45.10.Na	Geometrical and tensorial methods
% 	45.20.-d	Formalisms in classical mechanics
%		45.20.da	Forces and torques
%		45.20.dc	Rotational dynamics
%		45.40.Bb	Rotational kinematics
%
%-------------------------------------------------------------------------
% Parole chiave
\keywords{
%\textbf{cross product, pseudovector, N Dimensions, N-D, dimensional, wedge product, doublewedge}
\textbf{cross product, pseudovector, N Dimensions, dimensional, moment, N-D,  wedge product, doublewedge}
}
%------------------------------------------
%
%====================================================================
%
% Versione di REVTEX4.1
\maketitle
%
%==================================================================
% CONTROREVISIONE - DA USARE?
%==================================================================
%
\begin{comment}
% \fi
\IEEEpeerreviewmaketitle
\end{comment}
%
%===============================================================
%
%											*** MAINMATTER ***
%
%===============================================================
%
%\input{mainmatter_articolo_prl.tex}
%
%==================================================================
% MAINMATTER ARTICOLO PRL 2014
%==================================================================
%
%\input{schema_articolo_prl.tex}
%
%---------------------------------------------------------------------------------------------
%
\section{Introduction}
In this report we present a summary of a Master Thesis, published in Italian,
concerning the extension of cross product $\times$ in N Dimensions \cite{Gonano:tesi}. To indicate
that new operator we use the doublewedge $\doublewedge$ symbol, which resemble the Grassmann's wedge
product $\wedge$ and a notation for cross product commonly adopted in Italy (see also \cite{AmaldiCivita:mecc_razio}).
%CITARE LEVI CIVITA)
%
Here our task is to show the main drawbacks and difficulties of 3-D cross product and to
introduce a \textit{user-friendly} N-D notation, %an N-D notation of easy use,
suitable also for students.
%
%
%------------------------------------------------------------------------------------------------------
\subsection{Very brief historical notes}
\begin{figure}[htbp]
	\centering
		\includegraphics[width = 0.60\columnwidth]
		%{Figure/746px-Em_monopoles_svg.png}
	   	%{Figure/cross_wedge_product_rid.png}
%	   	{cross_wedge_product_rid.jpg}
	   	{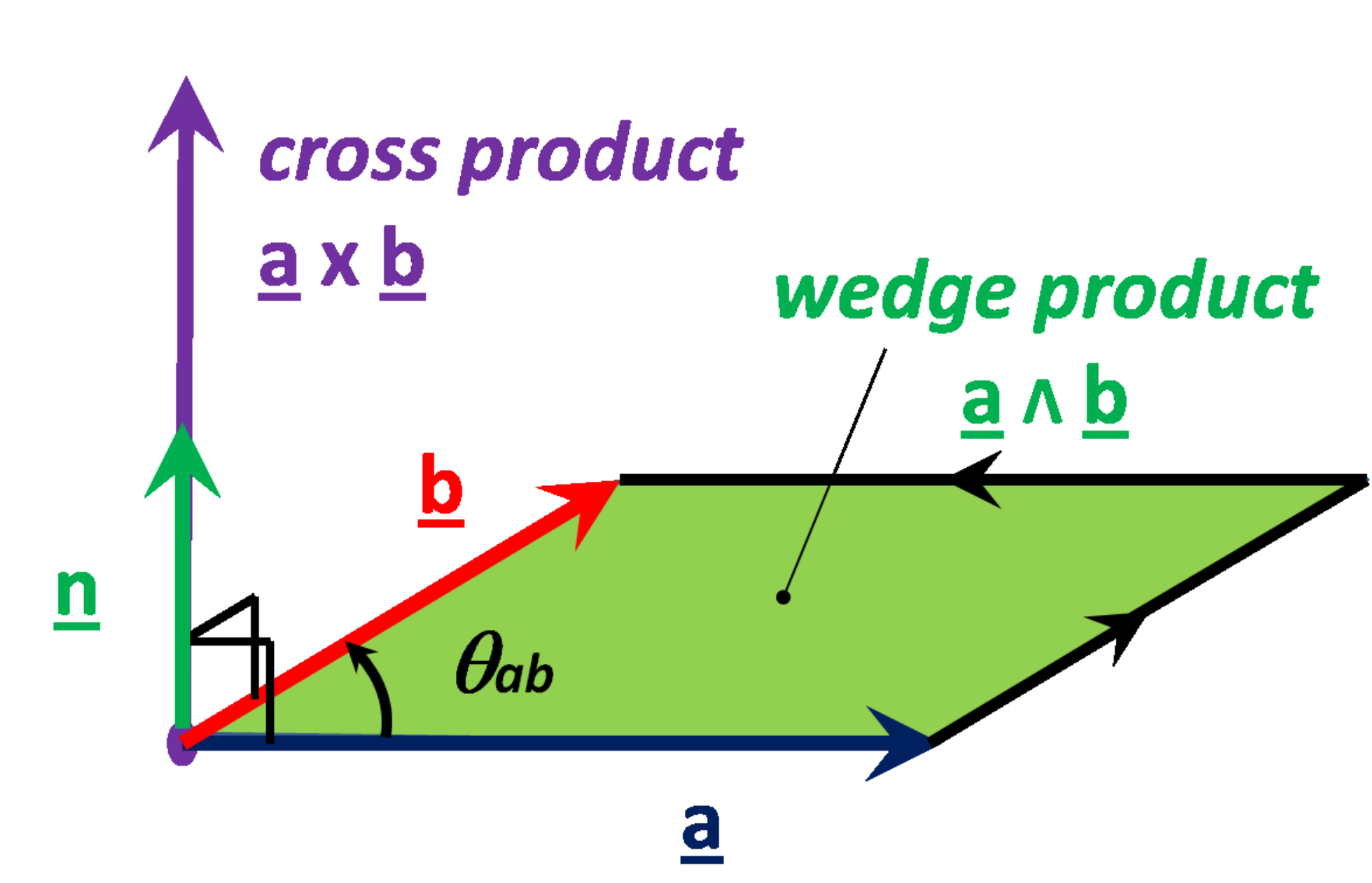}
	\caption{	\label{fig:cross_wedge_product}
	%\footnotesize
						Different interpretations of cross and wedge product}
\end{figure}

The history of cross product
is strictly related to that of vector calculus
\cite{Crowe:History_vector,
Crowe:History_vector_article}.
In 1773, Lagrange calculated the volume of a tetrahedra finding cross product via analysis, %CITARE LAGRANGE?,
but %the
\virga vectors\virgc haven't been invented yet.
In 1799, C. F. Gauss and C. Wessel represented complex numbers like arrows on a plane and %CITARE GAUSS E WESSEL ?.
in 1840 H.G. Grassmann introduced the \emph{exterior product} and a wedge $\wedge$ as its symbol. %CITARE GRASSMANN?.
That could be considered as the first cross product definition, but for Grassmann the operation's result
is not a \emph{vector}: though, it's an area or volume with an oriented boundary (Fig.\ref{fig:cross_wedge_product}).
In his External Algebra Grassmann also performs the first N-D extension of his operator $\wedge$, making it to act
on \underline{many} vectors at the same time%.
, e.g:
\begin{gather}
(\mathbf{a} \wedge \mathbf{b} \wedge \mathbf{d} ) \wedge \mathbf{c} = \mathbf{a} \wedge (\mathbf{b} \wedge \mathbf{c} \wedge \mathbf{d}) 
%\quad \text{associativity}
\end{gather}
The result of this operation is generally interpreted as the signed (hyper-)volume of a N-D parallelogram whose edges
are N vectors.
In 1843, W. R. Hamilton invented the quaternions to describe rotations in 3-D and in 1846 he adopted the terms
\emph{scalar} and \emph{vector} referring to real and imaginary parts of a quaternion.
The vector part of a product between quaternions with null real parts is equal to cross product.
%CITARE HAMILTON
%
In 1881-84, J.W. Gibbs wrote 
for his students the
\emph{Elements of Vector Analysis}\cite{Gibbs:Elements_of_vector_analysis},
where modern vector calculus is explained and
in 1901 his disciple E.B. Wilson published \emph{Vector Analysis}\cite{Wilson:Vector_analysis},
which had a large diffusion. In Gibbs's notation the
cross product is indicated with a $\times$ and it's considered a vector.
Shortly, from the end of the XIX century there were many different contributions to the development of
vector calculus, though interpretations and notations were not uniform. We can mention W.K. Clifford,
O. Heaviside, G. Peano, G. Ricci-Curbastro and T. Levi-Civita just to cite some who worked on that topic.
%
%----------------------------------------------------------------------------------------
%----------------------------------------------------------------------------------------
%
%Today
Nowadays cross $\times$ and wedge $\wedge$ products are well 
distinct
%distinguished %GRAZIE GIORGIO!
operators and 
employed
%used
in different fields, %contexts,
though they share similar algebraic properties.
%
%
%=============================================================================================
\section{3-D Cross product definition and uses}
The cross product is an operation between two vectors
$\vect{a}$ and $\vect{b}$ and in 3-D it is defined as $\vect{p} = \crossprod{a}{b}$ with:
\begin{gather}
%\prodvet{a}{b}
\vect{p}^T = \left[ a_2 b_3 - a_3 b_2,\, a_3 b_1 - a_1 b_3,\, a_1 b_2 - a_2 b_1 \right]
\end{gather}
The cross product frequently appears in Physics and Engineering, since it's used
for the calculus of moments, rotation axes, volumes, etc. 
\begin{align}
&\vect{M} = \crossprod{r}{F} \qquad \qquad \text{torque or moment of a force}\\
&\vect{\varphi}_{A \rightarrow B} = \varphi\, \frac{\crossprod{a}{b}}{\norma{\crossprod{a}{b}}}
%\varphi (\crossprod{a}{b})/\norma{\crossprod{a}{b}}
\quad \text{A-towards-B rotation-vector}\\
&V = (\crossprod{a}{b})^T\!\cdot\!\vect{c} \qquad \text{volume of parallelepiped $\vect{a},\vect{b},\vect{c}$}
\end{align}
Actually it's one of the most
widespread mathematical operator in Mechanics and it's suitable for many applications.
%
%In 3 Dimension $\vect{p}$ can be interpreted as a vector orthogonal to $\vect{a}$ and $\vect{b}$
%
\section{Limits and difficulties for 3-D cross product}
Though it is commonly used, the cross product presents some \virga oddities\virgc,
e.g., you need %it needs?
the concepts of \emph{clock-wise} sense and \emph{right-hand}.
%Moreover,
Furthermore,
%Particularly,
this operator is not always so easy to use:
%wrong signs are the most frequent mistake
the most frequent mistake is to confuse %wrong	%to make mistake with %GRAZIE GIORGIO
the signs
%Writing wrong signs is the most frequent mistake
%
($+$ or $-$ ?) and in pratice you have to \underline{memorize}
long identities like:
\begin{gather}
\vect{a} \times ( \crossprod{b}{c} )  =  ( \prodscal{c}{a} )\vect{b} -  ( \prodscal{b}{a} )\vect{c}\\
( \crossprod{a}{b} )^{T} \cdot ( \crossprod{c}{d} )  =
( \prodscal{a}{c} )( \prodscal{b}{d} ) - ( \prodscal{a}{d} )( \prodscal{b}{c} )
\end{gather}
Re-demostrate them every time is a long work, since it requires to explicit coordinates for each vector,
permutation of indices etc., and you risk confusion with letters and signs: mistake is in ambush.

%Moreover, there are more serious paradoxes, concerning cross product, which we are going to show.
Moreover, we are going to show some more serious paradoxes concerning cross product.
\subsection{3-D Rotation-vectors}
While sum and scalar product between vectors are operations easy to be extended in N-D,
the cross product is defined just in 3-D and it's often used to express \emph{rotation-vectors}.
Those kind of vectors can not be summed with the tip-tail rule, unlike common (polar) vector;
%
% INSERIRE QUI FIGURA DEI VETTORI ROTAZIONE?
%
in fact rotations don't sum because they don't commute:
\begin{gather}
%$\vect{\varphi}_{A \rightarrow B} +  \vect{\varphi}_{B \rightarrow C} \neq \vect{\varphi}_{A \rightarrow C}$
\vect{\varphi}_{A \rightarrow B} +  \vect{\varphi}_{B \rightarrow C} \neq \vect{\varphi}_{A \rightarrow C}
\end{gather}
%\\
Usually, those originated by cross product are called \emph{axial vectors} or \emph{pseudovectors}.
\subsection{Alice through the looking-glass}
If we place a set of \virga true\virgc vectors, like radii, velocities, forces etc., in front of the mirror
they are simply reflected, instead of moments and pseudovector in general.
%----------------------------------------------------------------------------------------
\begin{figure}[htbp]
	\centering
		\includegraphics[width = 0.80\columnwidth]
		%{Figure/cross_product_mirror.png}
%		{cross_product_mirror.jpg}
		{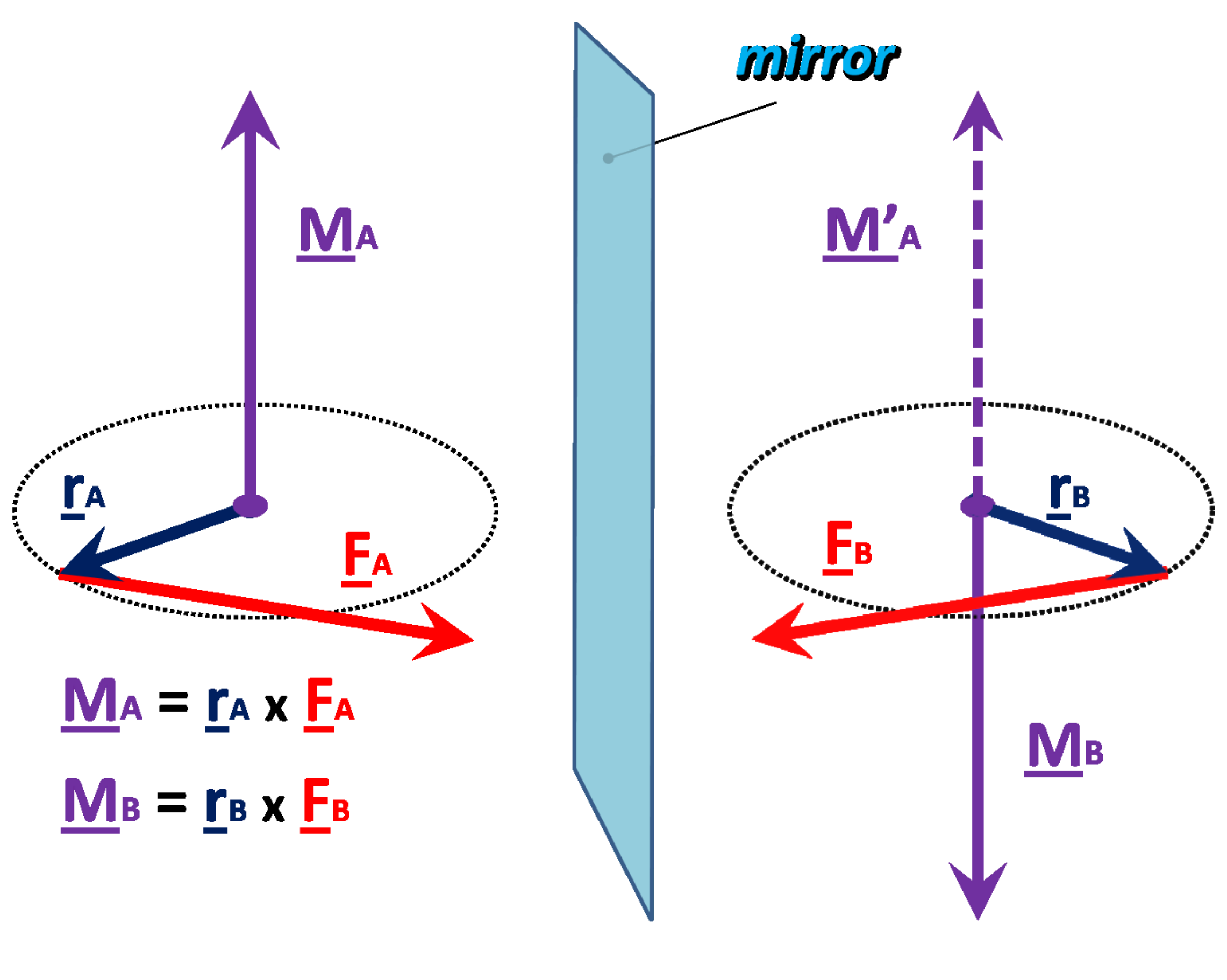}
	   %{Figure/vectors_mirror.png}
	\caption{\footnotesize Radius, force and moment's reflection}
	\label{fig:vectors_at_mirror}
\end{figure}
%----------------------------------------------------------------------------------------
%
In fact cross product doesn't respect reflection rules
and the specular image of a right hand is a left one and counterclock-wise looks clock-wise.
\subsection{Flatland - a 2-D world}
In \emph{Flatland}\cite{Abbott:Flatland} E. A. Abbott describes life and customs of people in a 2-D world:
in this universe
%In a 2-D world
vectors can be summed together and projected, 
areas are calculated, rotations are clock-wise or counterclock-wise, reflection is possible\dots
but cross product \emph{does not exist}; otherwise, 2-D inhabitants should have great fantasy to
imagine a $3^{rd}$ dimension to contain a vector orthogonal to their plane.
%
%----------------------------------------------------------------------------------------
\begin{figure}[htbp]
	\centering
		\includegraphics[width = 0.60\columnwidth]
		%{Figure/Flatland_(first_edition)_page_79.png}
%		{Flatland_(first_edition)_page_79.jpg}
		{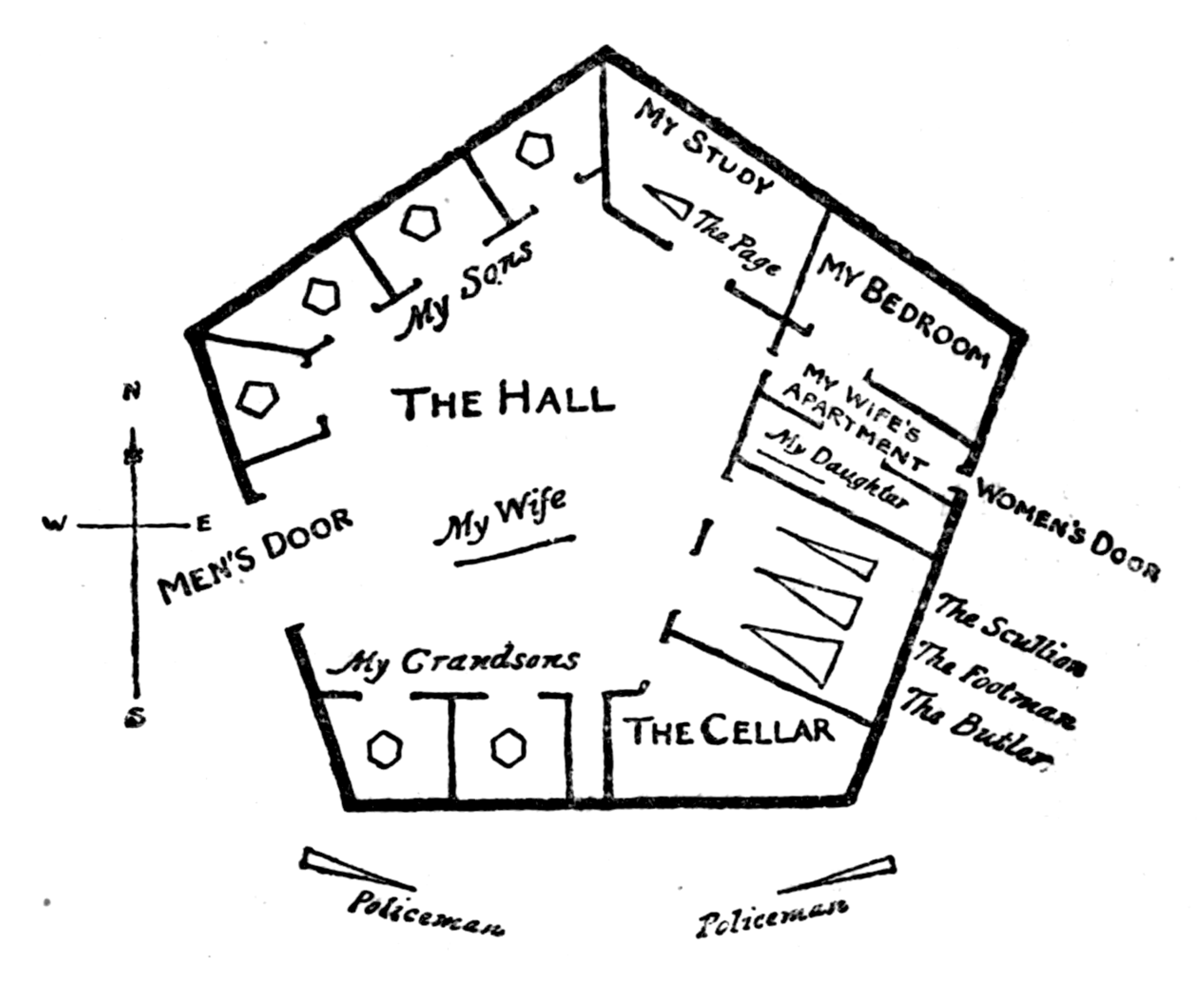}
	   %{Figure/vectors_mirror.png}
	\caption{\label{fig:flatland_hall}
	A picture from %Abbott's
	\emph{Flatland} - image in the public domain}
\end{figure}
%----------------------------------------------------------------------------------------
%
By the way, in 2-D a single scalar number is sufficient to describe a force's moment:
\begin{gather}
%\begin{equation}
M = M(\vect{r}, \vect{F}) = r_1 F_2 - r_2 F_1
%\end{equation}
\end{gather}
With such a definition, this operation respects \emph{all} algebraic properties of cross product,
but the result is a \emph{scalar}.
\subsection{4-D space}
In a 4-D space each vector has 4 components %elements
and in order to construct a cross product $\vect{p} = \crossprod{a}{b}$
we have to impose that $\vect{p}^T = \left[p_1, p_2, p_3, p_4\right]$ is perpendicular to vectors $\vect{a}$ and $\vect{b}$ 
and that its
magnitude is equal to the area between them:
\begin{gather}
%\vect{p}^T = \left[p_1, p_2, p_3, p_4\right]\\
\prodscal{p}{a} = 0; \,\,\, \prodscal{p}{b} = 0;
%\quad p^2 = a^2 b^2 - (\prodscal{a}{b})^2%quadrato aggiunto come correzione
\quad \norma{\vect{p}\,} = \norma{\,\crossprod{a}{b}}%quadrato aggiunto come correzione
%\begin{sistemalign}
%&\prodscal{p}{a} = 0 \\
%&\prodscal{p}{b} = 0 \\
%&\norma{p}^2 = \norma{a}^2 \norma{b}^2 - ( \prodscal{a}{b})^2
%p^2 = a^2 b^2 - (\prodscal{a}{b})^2
%\end{sistemalign}
%\norma{p}^2 = \norma{a}^2 \norma{b}^2 - \left( \prodscal{a}{b}\right)^2
\end{gather}
But these are just 3 equations, and we have 4 unknowns:
the problem has 1 degree of indetermination.
In fact, in 4-D there is an infinity of vectors $\vect{p}$ that satisfy these requirements: rotation axes are not unique!

So, cross product maybe exists just in 3-D, or it's not a vector.
%====================================================================================================
\section{N-D cross product}
As we have seen, in 3-D cross product can give some troubles. Now we desire to extend it in N
spatial Dimensions and want it %that %new operator
to satisfy some conditions:
\begin{itemize}
\item \textbf{Moment calculus}: the new operation should involve just 2 vectors %per
at time.

In fact, differently from the %Grassmann
exterior product $\wedge$, it must be %useful
of practical utility
in Physics
for calculating moments rather than volumes or determinants.
%for the calculus of moments rather than for that of volumes or determinants.
%
\item \textbf{Analogy}: the algebraic properties
of the new operator should be analogous to those
of the classic 3-D cross product.
\item \textbf{N-D validity}:
the new operation must be valid in \underline{every} positive integer number N of spatial Dimensions. %(N-D validity)
\item \textbf{User-friendly}: the N-D notation should be %of easy use, allowing simpler counts
general and of easy use, %handy PROPOSTO DA GIORGIO in alternativa a "`easy use"'
allowing simpler counts.
\end{itemize}
Moreover, we would like to solve some of the paradoxes previously mentioned, re-interpreting the
concept of cross product itself.
%
%\subsection{N-D extension for cross product}
%\subsection{Extending cross product in N-D}
\subsection{Definition of N-D cross product}
We notice that in Mechanics the angular velocity is sometime written like a pseudo-vector $\vect{\omega}$,
other times like a matrix $\doubleoverline{\Omega}$, and the latter %which
can be constructed %calculated
also in N-D.
For example, for two points $P$ and $Q$ on a rigid body we can write
the velocities $\vect{v}$ as:
\begin{align}
&\vect{v}_P  - \vect{v}_Q = \vect{\omega} \times (\vect{x}_P - \vect{x}_Q) 										\qquad \text{3-D notation}\\
&\vect{v}_P  - \vect{v}_Q =  \,\, \doubleoverline{\Omega}\cdot (\vect{x}_P - \vect{x}_Q)		\qquad \text{N-D notation}
\end{align} 
%
%VERSIONE PRECEDENTE
%In 3-D we can pass from one notation to the other through the well-know
%relation:
%
%-----------------------------------------------------------------------------------------------
In 3-D it's possible to pass from one notation to the other using the Levi-Civita $\varepsilon_{ijk}$ 
anti-symmetric
3-tensor:
%Analitically, \eqref{eq:from_omega_to_Omega} can be also expressed with the Levi-Civita tensor $\varepsilon_{ijk}$:
\begin{gather}
\omega_i  =   - \frac{1}{2}\sum_{j=1}^3 \sum_{k=1}^3 \left( \varepsilon_{ijk}\, \Omega_{jk}\right)\\
\Omega_{ij} = - \sum_{k=1}^3  \left( \varepsilon_{ijk}\, \omega_{k}\right)
%\omega_i  =   - \frac{1}{2}\sum_{j=1}^3 \sum_{k=1}^3 \left( \varepsilon_{ijk}\, \Omega_{jk}\right)
%\, \leftrightarrow \,
%\Omega_{ij} = - \sum_{k=1}^3  \left( \varepsilon_{ijk}\, \omega_{k}\right)
%
\end{gather}
However, using a tensor of rank 3 could be heavy for somebody, so we can write more simply:
\begin{gather} 
\overline{\overline{\Omega}} = [\omega \times] = %
\begin{bmatrix}
0    & -\omega_3 & \quad\!\omega_2 \\
\quad\!\omega_3  & 	0	 & -\omega_1 \\
-\omega_2 & \quad\!\omega_1 & 0
\end{bmatrix}
\label{eq:from_omega_to_Omega}
\end{gather}
%
%-----------------------------------------------------------------------------------------------
Is it possible a similar %analogous
reasoning %formulation
with moments?

%-----------------------------------------------------------------------------------------------
%Observing 
Let's observe
the $z$-component of a moment $M_z = r_x F_y - r_y F_x$: we notice that subscript \emph{z}
doesn't appear neither in the force nor in the radius, so $\vect{M}$, rather than \emph{\virga around z axis\virgc},
looks to be \emph{\virga from x to y\virgc}.
If we assemble the moment in a matrix form, we get:
\begin{gather}
\overline{\overline{M}} = [M \times] = %
\begin{bmatrix}
0    & -M_3 & \quad\!M_2 \\
\quad\!M_3  & 	0	 & -M_1 \\
-M_2 & \quad\!M_1 & 0
\end{bmatrix} \\%=
\overline{\overline{M}} =
\begin{bmatrix}
0    							 &  r_2 F_1 - r_1 F_2  &  r_3 F_1 - r_1 F_3 \\
r_1 F_2 - r_2 F_1  & 	0	 								&  r_3 F_2 - r_2 F_3  \\
r_1 F_3 - r_3 F_1  & r_2 F_3 - r_3 F_2  & 0
\end{bmatrix}
\end{gather}
It's straightforward to demonstrate that:
\begin{gather}
%\doubleoverline{M}_{ij}= r_j F_i - r_i F_j
M_{ij}=  F_i r_j- r_i F_j
\end{gather}
Since vectors $\vect{F}$ and $\vect{r}$ can have any dimension N, we define
the \textbf{N-D cross product} as the difference of dyads: %dyadics:
\begin{gather}
\doubleoverline{M} = \prodvett{r}{F} = \prodvettdiad{r}{F}
\end{gather}
It can be easily verified that the new operator respects 
all the required algebraic properties; just the result is no more a vector but an anti-symmetric
matrix or 2-tensor.
For full theory, see \cite{Gonano:tesi}.
%
%---------------------------------------------------------------------------------------------
%
%\subsection{\LaTeX command for N-D cross product}
%\subsection{\LaTeX command for doublewedge $\doublewedge$ symbol}
\subsection{\LaTeX command for the doublewedge symbol}
In order to distinguish the N-D cross product from the 3-D $\times$ and the wedge %exterior
$\wedge$ ones,
we introduced the new symbol $\doublewedge$, called \virga \textbf{doublewedge}\virgc.
In order to write the doublewedge %$\doublewedge$
in \LaTeX, you can create (or copy-paste)
a \emph{macro} %special command 
in the document
preamble:
\begin{verbatim}
\newcommand{\doublewedge}{\:\overset{\wedge}%
{\text{\scriptsize{$\wedge$}}}\;}
\end{verbatim}
%
%Then, to write the symbol, %operator,
%just digit \verb+\doublewedge+. %FRASE MIGLIORABILE
Then, to display the symbol, %operator,
just write \verb+\doublewedge+. 
%
%
%
%Per scrivere il doppio cuneo $\doublewedge$ in \LaTeX %
%-----------------------------------------------------------------------------------------
%è necessario caricare il pacchetto \verb+amssymb+, che arricchisce la scelta di 
%simboli matematici, scrivendo nel preambolo del documento: \verb+\usepackage{amssymb}+
%\\
% NON E' VERO! HO VERIFICATO CHE amssymb non è necessario!
%-----------------------------------------------------------------------------------------
%Sempre nel preambolo conviene creare l'apposito comando \verb+\doublewedge+:\\
%conviene creare nel preambolo l'apposito comando \verb+\doublewedge+:
%\\
%\\
%
%\verb+\newcommand{\doublewedge}{\:\overset{\wedge}%+
%\verb+{\text{\scriptsize{$\wedge$}}}\;}+
%-----------------------------------------------------------------------------------------
%
%\section{Algebraic properties and identities}
\subsection{Algebraic properties}
The N-D cross product 
or \emph{doublewedge product}
has many algebraic properties in common with the %of the
3-D one, as previously required.
\begin{itemize}
\item anti-commutativity: %anti-simmetry: %
\begin{gather}
\prodvett{a}{b} = -\, \prodvett{b}{a}
%= - \left[\prodvett{a}{b}\right]^T
\end{gather}
\item distributivity over addition:
\begin{gather}
\vect{a}\,\doublewedge \!\left(\vect{b} + \vect{c}\right) = \prodvett{a}{b} +  \prodvett{a}{c} 
\end{gather}
\item compatibility with scalar multiplication:
\begin{gather}
\left(\alpha\vect{a}\right)\doublewedge\left(\beta\vect{b}\right)=
\alpha\,\beta\, [\prodvett{a}{b}]  \quad \forall \,\alpha, \beta \in \mathbb{C}%\quad
\end{gather}
\end{itemize}
\begin{comment}
%-----------------------------------------------------------------------------------------
%VERSIONE PRECEDENTE, PIU' BREVE 
\begin{align}
&\prodvett{a}{b} = -\, \prodvett{b}{a} %\qquad
&\text{anticommutativity} \\
%
&\vect{a}\,\doublewedge \!\left(\vect{b} + \vect{c}\right) = \prodvett{a}{b} +  \prodvett{a}{c} %\quad
&\text{distributivity
%over sum
} \\
%
&\left(\alpha\vect{a}\right)\doublewedge\left(\beta\vect{b}\right)=
\alpha\,\beta\, [\prodvett{a}{b}]  %\quad
&\text{scalar multiplication}
%
\end{align}
\end{comment}
%-----------------------------------------------------------------------------------------
Differently from the cross and wedge products, the $\doublewedge$ operation cannot be repeated over itself,
since its \emph{inputs} are vectors and the \emph{output} is a matrix.
%
%support associativity
%
%=============================================================================================
%TABELLA CON LE IDENTITA' ALGBERICHE?
%
\subsection{Main algebraic identities}
In %the two-column
table \ref{tab:vectorial_identities} we report the main mathematical identities involving
cross product with both 3-D and N-D notations.
%INSERIRE QUI TABELLA?
%\input{tabella_algebric_id.tex}
%
%========================================================
% TABELLA IDENTITA' ALGEBRICHE
%========================================================
%
%In the two-column table \ref{tab:vectorial_identities} we report the main mathematical identities involving
%cross product with both 3-D and N-D notations.
%
%\begin{widetext}
% TABELLA
\begin{table*}[htbp]%The best place to locate the table environment is directly after its first reference in text
\caption{%
\label{tab:vectorial_identities}%
Main mathematical identities for cross product %with 3-D and N-D notations
}
\begin{ruledtabular}
% "`Ruled tabular"` è una tabella con le doppie righe all'inizio e alla fine che la delimitano.
%\begin{tabular}{lcdr}
% lcdr = left, centered, decimal, right
\begin{tabular}{ll}%{cc}%{rl}
%-------------------------------------
\textrm{3-D notation}&
%\textrm{Centered}&
%\multicolumn{1}{c}{\textrm{Decimal}}&
\textrm{N-D notation}\\
\colrule	%Linea orizzontale
%%--------------------------------------------------
%$\crossprod{a}{b}$  																& $\prodvett{a}{b}$  \\
$\vect{M} = \crossprod{r}{F}$  											& $\doubleoverline{M} = \prodvett{r}{F} = \prodvettdiad{r}{F}$  \\
%$\vect{M} = \left[r_2 F_3 - r_3 F_2,\, r_3 F_1 - r_1 F_3,\, r_1 F_2 - r_2 F_1\right]^T$
%--------------------------------------------------
$M_x = F_z r_y - r_z F_y$														& $M_{ij} = F_i r_j  - r_i F_j$ \\
%%--------------------------------------------------
$\left( \crossprod{r}{F} \right) \times \vect{c}
= \vect{F} \,(\prodscal{r}{c})  - \vect{r}\,(\prodscal{F}{c})
$
& $\left[ \prodvett{r}{F} \right] \cdot \vect{c}
= \vect{F} \,(\prodscal{r}{c})  - \vect{r}\,(\prodscal{F}{c})
$ \\
%%--------------------------------------------------
$\vect{M} \times \vect{c} = \left[M \times\right] \vect{c}$ 	& $\doubleoverline{M} \cdot \vect{c}$ \\
%
%%--------------------------------------------------
$V =(\crossprod{a}{b}) \cdot \vect{c}$ & $V =\left[ \prodvett{a}{b} \right]   \underset{321}{\cdot} \vect{c} $\\
$(\crossprod{a}{b}) \cdot \vect{c} = (\crossprod{b}{c}) \cdot \vect{a} =(\crossprod{c}{a}) \cdot \vect{b}
$
&
$\left[\prodvett{a}{b}\right]\underset{321}{\cdot}\vect{c} = \left[\prodvett{b}{c}\right]\underset{321}{\cdot}\vect{a} =
 \left[\prodvett{c}{a}\right]\underset{321}{\cdot}\vect{b} 
$\\
%
%%--------------------------------------------------
%FORMA ESPLICITA DEI VOLUMI 3-D - SCOMMENTARE?
%$ V =(a_1 b_2 c_3 + a_2 b_3 c_1  + a_3 b_1 c_2) -  (a_1 b_3 c_2 +  a_2 b_1 c_3 +a_3 b_2 c_1 )  $ &
%$V =(b_3 a_2  - a_3 b_2) c_1    + (b_1 a_3  - a_1 b_3) c_2      + (b_2 a_1 - a_2 b_1) c_3 $\\
%--------------------------------------------------
%LINEAR TRASFORMATION
$ ( \doubleoverline{L}\,\vect{a} ) \times ( \doubleoverline{L}\,\vect{b} ) =
\det(\doubleoverline{L}) %\left(\doubleoverline{L}^{-T}\, (\crossprod{a}{b}) \right) $ &
\left(L^{-T}\cdot (\crossprod{a}{b}) \right) $ &
$%\left( \doubleoverline{L}\,\vect{a} \right) \doublewedge \left( \doubleoverline{L}\,\vect{b} \right) =
( \doubleoverline{L}\,\vect{a} ) \doublewedge ( \doubleoverline{L}\,\vect{b} ) =
\doubleoverline{L} \,\left[ \prodvett{a}{b} \right] \,\doubleoverline{L}^T$\\
%
%--------------------------------------------------
$(\crossprod{a}{b}) \cdot (\crossprod{c}{d})
%=( \prodscal{a}{c} )( \prodscal{b}{d} ) - ( \prodscal{a}{d} )( \prodscal{b}{c} )
=( \vect{a}\cdot\vect{c} )( \vect{b}\cdot\vect{d} ) - ( \vect{a}\cdot\vect{d} )(\vect{b}\cdot\vect{c} )
$ & 
$ \frac{1}{2} \left[ \prodvett{a}{b} \right] :  \left[ \prodvett{c}{d} \right]
=( \prodscal{a}{c} )( \prodscal{b}{d} ) - ( \prodscal{a}{d} )( \prodscal{b}{c} )
$ \\
%--------------------------------------------------
$Pow = \vect{M} \cdot \vect{\omega}$ & 
$Pow = \frac{1}{2}\doubleoverline{M} : \doubleoverline{\omega}$\\
%
%\vect{p}^T = \left[
%a_2 b_3 - a_3 b_2,\, a_3 b_1 - a_1 b_3,\, a_1 b_2 - a_2 b_1
%\right]
%
%-------------------------------------
\end{tabular}
\end{ruledtabular}
\end{table*}
\section{Applications and consequences}
In this section we bring %report
some sparse examples regarding the application of $\doublewedge$ product in different
contexts. %frameworks? fields? sectors?
For details see\cite{Gonano:tesi}.
%
%------------------------------------------------------------------------------------------------------
\subsection{Perpedicular component of a vector}
The perpendicular component $\vect{F}_{\perp}$ of a vector $\vect{F}$ on an other $\vect{r}$ can be calculated as:
\begin{gather}
\vect{F}_{\perp} = \vect{F}  -  \frac{1}{r^2}( \prodscal{r}{F}) \vect{r}
\quad \Longrightarrow \quad \vect{r}^T \!\cdot \! \vect{F}_{\perp} =0
\end{gather}
The same equation can be re-written as:
\begin{gather}
\vect{F}_{\perp} = \frac{1}{r^2}\left( \proddiad{F}{r}  -  \proddiad{r}{F}\right)\cdot\vect{r} =
\frac{1}{r^2} \,[ \prodvett{r}{F}] \cdot \vect{r}
\end{gather}
This result is a particular case of the identity:
\begin{align}
[\prodvett{r}{F}] \cdot \vect{c} &= \left( \proddiad{F}{r}  -  \proddiad{r}{F}\right)\cdot\vect{c} \\
[\prodvett{r}{F}] \cdot \vect{c} &= \vect{F} \,(\prodscal{r}{c})  - \vect{r}\,(\prodscal{F}{c})
\end{align}
Let's notice that we derived %have drawn
it in 2 rows. The 3-D equivalent identity is:
\begin{gather}
\left( \crossprod{r}{F} \right) \times \vect{c} = \vect{F} \,(\vect{r}\cdot\vect{c})  - \vect{r}\,(\vect{F}\cdot\vect{c})
\end{gather}
but to demonstrate it with 3-D formalism %notation
it's a longer task (try to believe).
%
%------------------------------------------------------------------------------------------------------
%\subsection{Mechanic applications}
\subsection{Angular momenta and inertia matrices}
Given a body defined on a lagrangian domain $\Omega_x$, its angular momentum $\doubleoverline{L}_0$ with respect to
a pole $\vect{x}_0$ is:
\begin{gather}
\doubleoverline{L}_0 = \int_{\Omega_x} (\vect{x} - \vect{x}_0 ) \doublewedge (\rho\vect{v})\, d\Omega_x
\end{gather}
where $\rho$ and $\vect{v}$ are the mass density and velocity respectevely.
In 3-D, for a rigid body holds: %stands:
\begin{gather}
\vect{L}_0 = m\,(\vect{x}_G - \vect{x}_0) \times \vect{v}_0 + \doubleoverline{\mathcal{I}}_0 \,\vect{\omega}
\label{eq:angular_moment_rigid_body_3D}
\end{gather}
where $\vect{x}_G$ is the center of gravity and $\doubleoverline{\mathcal{I}}_0$ is the 3-D inertia matrix, defined as:
\begin{equation}
\doubleoverline{\mathcal{I}}_0 =%\triangleq
\int_{\Omega_x} \left( 
\rho 
\begin{bmatrix}
  y^2 + z^2 &  - x y 				& 	- z x \\
 -x y				&   z^2 + x^2	  & 	-  y z  \\
 -z x 			&  - y z 				&		 x^2 + y^2
\end{bmatrix}  
\right)
\; d\Omega_x 
\label{eq:3_D_inertia_matrix}
\end{equation}
where $[x; y; z] = \vec{x} - \vec{x}_0$.
Let's notice that in \eqref{eq:3_D_inertia_matrix} indices are misleading, in fact:
\begin{gather}
\mathcal{I}_{xx}= \int_{\Omega_x} \rho\, (y^2 + z^2)\; d\Omega_x
\neq \int_{\Omega_x} \rho \,x^2\; d\Omega_x
\end{gather}
With N-D notation, instead, the inertia matrix $\doubleoverline{I}_0$ is compactly defined as:
\begin{equation}
\doubleoverline{I}_0 \triangleq
\int_{\Omega_x} \rho\, \proddiad{\Delta{x}_0}{\Delta{x}_0}\; d\Omega_x
\qquad \text{with: }\Delta\vect{x}_0 = \vect{x} - \vect{x}_0
\end{equation}
Let's notice that the N-D inertia matrix 
$\doubleoverline{I}_0$
is conceptually similar to the matrix of covariances $\sigma_{ij}^2$
used in Statistics.

With N-D notation the Eq.\eqref{eq:angular_moment_rigid_body_3D} will look:
\begin{gather}
\doubleoverline{L}_0 = m\,(\vect{x}_G - \vect{x}_0) \doublewedge \vect{v}_0 +
\doubleoverline{I}_0 \,\doubleoverline{\Omega} -\left( \doubleoverline{I}_0 \,\doubleoverline{\Omega}\right)^T
\end{gather}
%
%------------------------------------------------------------------------------------------------------
\subsection{Volume calculus: the 3-indices product}
In 3-D, the signed volume $V$ of a parallelepiped whose edges are vectors $\vect{a}, \vect{b}, \vect{c}$
can be calculated as:
\begin{gather}
V = (\crossprod{a}{b}) \cdot \vect{c} = (\crossprod{b}{c}) \cdot \vect{a} =(\crossprod{c}{a}) \cdot \vect{b} \\
\begin{split}
V =\quad &(b_3 a_2  - a_3 b_2) c_1    + (b_1 a_3  - a_1 b_3) c_2  \\
   + &(b_2 a_1 - a_2 b_1) c_3
\end{split}
\end{gather}
In N-D for a hyper-parallelepiped with edges $\vect{v}_1, \vect{v}_2,\cdots \vect{v}_N$
the signed hyper-volume is:
\begin{gather}
V = \det \left| \vect{v}_1, \vect{v}_2,\cdots \vect{v}_N \right|
\end{gather}
However, if we want to determine a 3-D volume in an N-D space it's convenient to define the \emph{3-indices product}:
\begin{eqnarray}
\doubleoverline{A} \underset{ijk}{\cdot} \vect{c} = A_{ij} \,c_k + A_{jk} \,c_i + A_{ki} \,c_j %\quad
%\forall i,j,k \in \left{1,2,\cdots, N\right}\\
\\
\doubleoverline{A} \underset{ijk}{\cdot} \vect{c} =
%\doubleoverline{A} \cdot_{ijk} \vect{c} =
\doubleoverline{A} \underset{jki}{\cdot}\vect{c} =
%\doubleoverline{A} \overset{\cdot}{\text{\scriptsize{jki}}}\vect{c} =
\doubleoverline{A} \underset{kij}{\cdot} \vect{c}
\end{eqnarray}
where $i,j,k \in\left\{1,2,\cdots, N\right\}$ are arbitrary indices.
It's quite straightforward to verify %demonstrate
that in 3-D holds:
\begin{equation}
V = (\crossprod{a}{b}) \cdot \vect{c} = \left[\prodvett{a}{b}\right]\underset{321}{\cdot}\vect{c} 
\end{equation}
Anyway we remember that the $\doublewedge$ operator was conceived for the calculus of moments rather than volumes.
%
%------------------------------------------------------------------------------------------------------
\subsection{Power calculus: the matrix contraction}
In Mechanics the power $P$ transferred to a rotating body is the scalar product of its angular velocity $\vect{\omega}$
and the applied torque $\vect{M}$
\begin{gather}
P = \vect{M} \cdot \vect{\omega}
\end{gather}
Since both $\vect{\omega}$ and $\vect{M}$ are pseudovectors, with the N-D formalism the power will be calculated by the
contraction $:$ of matrices $\doubleoverline{M}$ and $\doubleoverline{\omega}$
\begin{equation}
P = \frac{1}{2}\,\doubleoverline{M} : \doubleoverline{\omega} =
\frac{1}{2}\,\sum_{i=1}^N \sum_{j=1}^N M_{ij} \, \Omega_{ij}
\end{equation}
The basic idea is quite similar to the tensor contraction adopted in Relativity.
%
%
%------------------------------------------------------------------------------------------------------
\subsection{3-D and N-D curl}
Curl is the differential operator analogous to cross product, and in 3-D it suffers for the same
problems, since it generates pseudovectors.
\begin{gather}
%\!\vectnabla \!\wedge \!\vect{v} =%\!
\!\vectnabla \!\times \!\vect{v} =%\!
\left[
\gradmat{v_3}{x_2} \! -\! \gradmat{v_2}{x_3}; %\,
\gradmat{v_1}{x_3} \! -\! \gradmat{v_3}{x_1}; %\,
\gradmat{v_2}{x_1} \! -\! \gradmat{v_1}{x_2}
\right]
\end{gather}
The extension in N-D is instantaneous:
\begin{gather}
\rotorr{v} =\rotorrdiad{v}{x}\\
[\,\rotorr{v}\,]_{ij} =\gradmat{v_i}{x_j}  - \gradmat{v_j}{x_i} = v_{i/j} - v_{j/i}%aggiunta rispetto all'articolo sui monopoli
\end{gather}
Even in this case it can be verified that the new operator satisfies \emph{all}
the
required differential properties.
%
%=============================================================================================
%\section{Maxwell's Equations in N-D}
\subsection{The magnetic field B is not a vector}
The magnetic field $\bm{B}$ is often involved with cross product and curl: is it a \virga true\virgc vector?
Look at Faraday's law and Lorentz force equations in 3-D:
\begin{gather}
\vectnabla \times \vect{E} = - \gradmat{\vect{B}}{t} \quad \qquad \vect{F}_B = - Q_e \, \vect{B} \times \vect{v}
\label{eq:Faraday_Lorentz_eqs_3D}
\end{gather}
We know, \emph{from definition}, that $\vect{E}$, $\vect{F}_B$ and $\vect{v}$ are true vectors and,
using N-D notation, \eqref{eq:Faraday_Lorentz_eqs_3D} will look:
\begin{gather}
\rotorr{E} = -  \gradmat{\doubleoverline{B}}{t} \qquad \qquad \vect{F}_B = -  Q_e\,\doubleoverline{B} \,\vect{v}
\label{eq:Faraday_Lorentz_eqs_ND}
\end{gather}
Thus, the magnetic field $\bm{B}$ is a not a vector, but a \emph{pseudovector}, and, in a wider N-D view, it is a 
\emph{matrix} or \emph{2-tensor}.
The use of $\bm{B}$-tensor is not new, but it seems not to be always understood:
for further details see \cite{Gonano:magnetic_monopoles_ND}
%
%-------------------------------------------------------
\begin{comment}
%
\subsection{Divergence and curl for B field}
In 3-D the divergence of a pseudo-vector is:
\begin{align}
\divergv{w} &= %\frac{\partial w_x}{ \partial x} + \frac{\partial w_y}{ \partial y} +\frac{\partial w_z}{ \partial z}
\gradmat{w_x}{x}+ \gradmat{w_y}{y}+\gradmat{w_z}{z}\\
\divergv{w} &=  w_{32/1} + w_{21/3} + w_{13/2}
\end{align}
In a similar way, in N-D notation stands:
%
\begin{gather}
\nabla_{ijk}\cdot\doubleoverline{w} = w_{ij/k} + w_{jk/i} + w_{ki/j} 
\end{gather}
So $\divergv{w} = \nabla_{321}\cdot\doubleoverline{w}$, 
where $\nabla_{ijk}$ is called \emph{3-indices divergence}\cite{Gonano:tesi}.
It can be easily demonstrated that  3-D curl of a pseudo-vector  correspond to:
\begin{gather}
\rotor{B} = -[B \wedge] \vectnabla = - \doubleoverline{B} \cdot \vectnabla
\end{gather}
In general, curl of a pseudo-vector is a true vector.
%
\end{comment}
%
%
%=============================================================================================
% CONCLUSIONI
%=============================================================================================
%
\section{Conclusion}
With the usual 3-D notation the cross product exhibits many limits and difficulties, since it produces
pseudovectors. In order to simplify calculations we defined the N-D cross product and introduced
the $\doublewedge$ symbol, solving some paradoxes and showing that moments are actually better
described by matrices rather than by vectors.
In this paper we reported just a summary of a more complete work \cite{Gonano:tesi}
which also includes the N-D curl extension.
We underline that the use of 2-tensors instead of pseudovectors
\footnote{e.g. for magnetic field \textbf{B} and angular momentum \textbf{L}}
is not a completely new idea,
but it seems not to be so widespread or understood, even in Relativity and Quantum Mechanics.

The N-D notation for cross product was explicity
conceived to
%designed to %help the students
%help the students
help students with counts
%simplify the work of students,
and we are %quite
confident that it will be 
%a valid mathematical
a practical
tool also in classic Mechanics and Geometry.
%
%
%We analysed the limits and difficulties of the classic 3-D cross product and showed that
%it produces pseudo vectors. In order to solve the related paradoxes we 
%
% while with a N-D notation the generate matrices. We defined the cross product in troduce the $\doublewedge$ symbol and , 
%=============================================================================================
%
%\subsection{Abbreviations and Acronyms}
%\subsection{Units and Numbers}
%
%\subsection{References}
%\section*{ACKNOWLEDGEMENT}
%
%------------------------------------------------------------------
% RICONOSCIMENTI
%------------------------------------------------------------------
%
%\begin{acknowledgments}
\section*{Acknowledgments} %<---PER PRL I RINGRAZIAMENTI VANNO LASCIATI COME UN UNICO
%																PARAGRAFO DOPO LE CONCLUSIONI, SENZA TITOLO.
% Acknowledgments (if any) should appear as a separate non-numbered section
% before the list of references.
%
%We are grateful to Prof. Antonella Abb\'a for her kind advice and support.
%
%We thanks Prof.sa Sonia Leva, Riccardo Albi, Giorgio Fumagalli and 
We thanks Prof. Antonella Abb\'a, Prof. Sonia Leva,
Riccardo Albi, Giorgio Fumagalli, Andrea Gatti, Pietro Giuri and 
Alessandro Niccolai
\begin{comment}
%Emanuele Ogliari?
%Francesco Gant?
%Daniele Gatti?
%
i Cavalieri delle Zodiaco, Iron Man, l'Uomo Conundrum,
il tacchino di Russell, il cigno nero \dots
\end{comment}
%Prof.sa Chiara Bisagni?
%
for their careful reviews %revisions
and Prof. Marco Mussetta for his precious help and support. % with \LaTeX.
We also would like to signal some authors who have independently come %arrived
to conclusions
analogous to ours in different ways:
%
%\cite{McDavid_McMullen:Maxwell_ND}
%\cite{Guio:Levi_Civita_tensor}
\cite{McDavid_McMullen:Maxwell_ND,
			Guio:Levi_Civita_tensor}
\normalsize
%
% Versione con BibTex
%\input{bibliografia_Bibtex_prl.tex}
%
%============================================================
% BIBLIOGRAFIA BIBTEX
%============================================================
% Per info vedi "`ArteLaTeX"' di Lorenzo Pantieri, pag 169/248 e seguenti.
%------------------------------------------------------------
% References section
%
%
% Use the following option to include external BibTeX files:
%
%----------------------------
% Stile bibliografico
%----------------------------
% Lo "`Stile"' della citazione varia da rivista a rivista.
% In genere è definito in un file con ESTENSIONE .bst
%
%\bibliographystyle{IEEEtran}
%\bibliographystyle{plain} %classico		
%\bibliographystyle{plain_ita} %classico italiano
%
% Stile per riviste APS:
% You should use BibTeX and apsrev.bst for references
% Choosing a journal automatically selects the correct APS
% BibTeX style file (bst file), so only uncomment the line
% below if necessary.
\bibliographystyle{apsrev4-1} %APS Physic Review
%\bibliographystyle{apsrmp4-1} %APS Review of Modern Physics
%----------------------------------------------------------
%
%
%%Caricamento dei dati bibliografici
%\bibliography{biblio_Bibtex_prl_201406} %nome del file .bib, senza estensione
%
%Versione con database multipli:
%\bibliography{file_bib1, file_bib2,... etc.} nomi dei file .bib, senza estensione
%
%Opzione per riportare nella bibliografia anche le opere non citate esplicitamente
\nocite{*}
%
% N.B.: in genere questo comando viene posto prima di %\bibliographystyle{...}%
% VERSIONE IEEE (con BibTeX o \bibitem)
%\bibliographystyle{IEEEtran}

\begin{thebibliography}{18}%
%
\makeatletter
\providecommand \@ifxundefined [1]{%
 \@ifx{#1\undefined}
}%
\providecommand \@ifnum [1]{%
 \ifnum #1\expandafter \@firstoftwo
 \else \expandafter \@secondoftwo
 \fi
}%
\providecommand \@ifx [1]{%
 \ifx #1\expandafter \@firstoftwo
 \else \expandafter \@secondoftwo
 \fi
}%
\providecommand \natexlab [1]{#1}%
\providecommand \enquote  [1]{``#1''}%
\providecommand \bibnamefont  [1]{#1}%
\providecommand \bibfnamefont [1]{#1}%
\providecommand \citenamefont [1]{#1}%
\providecommand \href@noop [0]{\@secondoftwo}%
\providecommand \href [0]{\begingroup \@sanitize@url \@href}%
\providecommand \@href[1]{\@@startlink{#1}\@@href}%
\providecommand \@@href[1]{\endgroup#1\@@endlink}%
\providecommand \@sanitize@url [0]{\catcode `\\12\catcode `\$12\catcode
  `\&12\catcode `\#12\catcode `\^12\catcode `\_12\catcode `\%12\relax}%
\providecommand \@@startlink[1]{}%
\providecommand \@@endlink[0]{}%
\providecommand \url  [0]{\begingroup\@sanitize@url \@url }%
\providecommand \@url [1]{\endgroup\@href {#1}{\urlprefix }}%
\providecommand \urlprefix  [0]{URL }%
\providecommand \Eprint [0]{\href }%
\providecommand \doibase [0]{http://dx.doi.org/}%
\providecommand \selectlanguage [0]{\@gobble}%
\providecommand \bibinfo  [0]{\@secondoftwo}%
\providecommand \bibfield  [0]{\@secondoftwo}%
\providecommand \translation [1]{[#1]}%
\providecommand \BibitemOpen [0]{}%
\providecommand \bibitemStop [0]{}%
\providecommand \bibitemNoStop [0]{.\EOS\space}%
\providecommand \EOS [0]{\spacefactor3000\relax}%
\providecommand \BibitemShut  [1]{\csname bibitem#1\endcsname}%
\let\auto@bib@innerbib\@empty
%</preamble>
\bibitem [{\citenamefont {Gonano}(2011)}]{Gonano:tesi}%
  \BibitemOpen
  \bibfield  {author} {\bibinfo {author} {\bibfnamefont {C.~A.}\ \bibnamefont
  {Gonano}},\ }\emph {\bibinfo {title} {{Estensione in N-D di prodotto vettore
  e rotore e loro applicazioni}}},\ \href
  {http://hdl.handle.net/10589/34061?locale=en} {Master's thesis},\ \bibinfo
  {school} {Politecnico di Milano} (\bibinfo {year} {2011})\BibitemShut
  {NoStop}%
\bibitem [{\citenamefont {Levi-Civita}\ and\ \citenamefont
  {Amaldi}(1949)}]{AmaldiCivita:mecc_razio}%
  \BibitemOpen
  \bibfield  {author} {\bibinfo {author} {\bibfnamefont {T.}~\bibnamefont
  {Levi-Civita}}\ and\ \bibinfo {author} {\bibfnamefont {U.}~\bibnamefont
  {Amaldi}},\ }\href@noop {} {\emph {\bibinfo {title} {{Lezioni di meccanica
  razionale}}}},\ Vol.~\bibinfo {volume} {I}\ (\bibinfo  {publisher}
  {Zanichelli editore Bologna},\ \bibinfo {year} {1949})\BibitemShut {NoStop}%
\bibitem [{\citenamefont {Crowe}(1967)}]{Crowe:History_vector}%
  \BibitemOpen
  \bibfield  {author} {\bibinfo {author} {\bibfnamefont {M.~J.}\ \bibnamefont
  {Crowe}},\ }\href@noop {} {\emph {\bibinfo {title} {{A History of Vector
  Analysis: The Evolution of the Idea of a Vectorial System}}}}\ (\bibinfo
  {publisher} {University of Notre Dame press},\ \bibinfo {year}
  {1967})\BibitemShut {NoStop}%
\bibitem [{\citenamefont {Crowe}(2002)}]{Crowe:History_vector_article}%
  \BibitemOpen
  \bibfield  {author} {\bibinfo {author} {\bibfnamefont {M.~J.}\ \bibnamefont
  {Crowe}},\ }\href
  {http://www.math.ucdavis.edu/~temple/MAT21D/SUPPLEMENTARY-ARTICLES/Crowe_His%
tory-of-Vectors.pdf} {\enquote {\bibinfo {title} {{A History of Vector
  Analysis}},}\ } (\bibinfo {year} {2002}),\ \bibinfo {note} {talk at
  University of Louisville}\BibitemShut {NoStop}%
\bibitem [{\citenamefont {Gibbs}(1884)}]{Gibbs:Elements_of_vector_analysis}%
  \BibitemOpen
  \bibfield  {author} {\bibinfo {author} {\bibfnamefont {J.~W.}\ \bibnamefont
  {Gibbs}},\ }\href
  {http://www.archive.org/stream/elementsvectora00gibbgoog\string#page/n4/mode/2up}
  {\enquote {\bibinfo {title} {{Elements of Vector Analysis - Arranged for the
  Use of Students of Physics}},}\ } (\bibinfo {year} {1881-1884}),\ \bibinfo
  {note} {note for students, privately printed}\BibitemShut {NoStop}%
\bibitem [{\citenamefont {Wilson}(1901)}]{Wilson:Vector_analysis}%
  \BibitemOpen
  \bibfield  {author} {\bibinfo {author} {\bibfnamefont {E.~B.}\ \bibnamefont
  {Wilson}},\ }\href
  {http://www.archive.org/details/vectoranalysiste00gibbiala} {\emph {\bibinfo
  {title} {{Vector analysis - A text-book for the use of students of
  mathematics and physics}}}}\ (\bibinfo  {publisher} {Yale Bicentennial
  publication},\ \bibinfo {year} {1901})\BibitemShut {NoStop}%
\bibitem [{\citenamefont {Abbott}(1884)}]{Abbott:Flatland}%
  \BibitemOpen
  \bibfield  {author} {\bibinfo {author} {\bibfnamefont {E.~A.}\ \bibnamefont
  {Abbott}},\ }\href@noop {} {\emph {\bibinfo {title} {{Flatland: A Romance of
  Many Dimensions}}}}\ (\bibinfo  {publisher} {Seely \& Co.},\ \bibinfo {year}
  {1884})\BibitemShut {NoStop}%
\bibitem [{\citenamefont {Gonano}\ and\ \citenamefont
  {Zich}(2013)}]{Gonano:magnetic_monopoles_ND}%
  \BibitemOpen
  \bibfield  {author} {\bibinfo {author} {\bibfnamefont {C.~A.}\ \bibnamefont
  {Gonano}}\ and\ \bibinfo {author} {\bibfnamefont {R.~E.}\ \bibnamefont
  {Zich}},\ }
  %
  {\enquote {\bibinfo {title} {{Magnetic monopoles and Maxwell's equations in N Dimensions}},}\ }
  %
  %
   in\  \href {\doibase 10.1109/ICEAA.2013.6632510} {\emph {\bibinfo
  {booktitle} {Electromagnetics in Advanced Applications (ICEAA), 2013
  International Conference on}}}\ (\bibinfo {year} {2013})\ pp.\ \bibinfo
  {pages} {1544--1547}\BibitemShut {NoStop}%
\bibitem [{Note1()}]{Note1}%
  \BibitemOpen
  \bibinfo {note} {E.g. for magnetic field \protect \textbf {B} and angular
  momentum \protect \textbf {L}}\BibitemShut {NoStop}%
\bibitem [{\citenamefont {McDavid}\ and\ \citenamefont
  {McMullen}(2006)}]{McDavid_McMullen:Maxwell_ND}%
  \BibitemOpen
  \bibfield  {author} {\bibinfo {author} {\bibfnamefont {A.}~\bibnamefont
  {McDavid}}\ and\ \bibinfo {author} {\bibfnamefont {C.}~\bibnamefont
  {McMullen}},\ }\href {http://arxiv.org/ftp/hep-ph/papers/0609/0609260.pdf}
  {\enquote {\bibinfo {title} {{Generalizing Cross Products and Maxwell's
  Equations to Universal Extra Dimensions}},}\ } (\bibinfo {year} {2006}),\
  \Eprint {http://arxiv.org/abs/hep-ph/0609260} {arXiv:hep-ph/0609260 [hep-ph]}
  \BibitemShut {NoStop}%
%%CITATION = HEP-PH/0609260;%%
\bibitem [{\citenamefont {Guio}(2011)}]{Guio:Levi_Civita_tensor}%
  \BibitemOpen
  \bibfield  {author} {\bibinfo {author} {\bibfnamefont {P.}~\bibnamefont
  {Guio}},\ }\href
  {http://www.homepages.ucl.ac.uk/~ucappgu/seminars/levi-civita.pdf} {\enquote
  {\bibinfo {title} {{Levi-Civita symbol and cross product vector/tensor}},}\ }
  (\bibinfo {year} {2011}),\ \bibinfo {note} {original note developed for a
  course on Physics of Astrophysics}\BibitemShut {NoStop}%
  %
  %
  %
\bibitem [{\citenamefont {Gray}(1969)}]{Gray:vector_cross_prod}%
  \BibitemOpen
  \bibfield  {author} {\bibinfo {author} {\bibfnamefont {A.}~\bibnamefont
  {Gray}},\ }
%
  %
  {\enquote {\bibinfo {title} {{Vector Cross Products on Manifolds}},}\ }
  %
 % 
  \href {http://www.jstor.org/stable/1995115} {\bibfield  {journal}
  {\bibinfo  {journal} {Transactions of the American Mathematical Society}\
  }\textbf {\bibinfo {volume} {141}},\ \bibinfo {pages} {pp. 465-504} (\bibinfo
  {year} {1969})}\BibitemShut {NoStop}%
\bibitem [{\citenamefont {Hage-Hassan}(2006)}]{Hassan:inertia_tensor_ND}%
  \BibitemOpen
  \bibfield  {author} {\bibinfo {author} {\bibfnamefont {M.}~\bibnamefont
  {Hage-Hassan}},\ }\href {http://arxiv.org/abs/math-ph/0604051v1} {\enquote
  {\bibinfo {title} {{Inertia tensor and cross product In n-dimensions
  space}},}\ } (\bibinfo {year} {2006}),\ \Eprint
  {http://arxiv.org/abs/math-ph/0604051} {arXiv:math-ph/0604051 [math-ph]}
  \BibitemShut {NoStop}%
%%CITATION = MATH-PH/0604051;%%
\bibitem [{\citenamefont {Manarini}(1939)}]{Manarini:estensione_ND_1939}%
  \BibitemOpen
  \bibfield  {author} {\bibinfo {author} {\bibfnamefont {M.}~\bibnamefont
  {Manarini}},\ }
 % 
  {\enquote {\bibinfo {title} {{Estensione della formula del doppio prodotto vettoriale agli spazi a pi\'u di tre dimensioni. Una formula di calcolo integrale ed un teorema della divergenza per i bivettori}},}\ }
%  
  %
  \href {http://www.numdam.org/item?id=RSMUP_1939__10__1_0}
  {\bibfield  {journal} {\bibinfo  {journal} {Rend. Semin. Mat. Univ. Padova}\
  }\textbf {\bibinfo {volume} {10}},\ \bibinfo {pages} {1-20} (\bibinfo {year}
  {1939})}\BibitemShut {NoStop}%
  %
  %
%
\bibitem [{\citenamefont {Morando}\ and\ \citenamefont
  {Leva}(1998)}]{Morando_Leva:note_vect}%
  \BibitemOpen
  \bibfield  {author} {\bibinfo {author} {\bibfnamefont {A.~P.}\ \bibnamefont
  {Morando}}\ and\ \bibinfo {author} {\bibfnamefont {S.}~\bibnamefont {Leva}},\
  }\href@noop {} {\emph {\bibinfo {title} {{Note di teoria dei Campi
  Vettoriali}}}}\ (\bibinfo  {publisher} {{Esculapio, Bologna}},\ \bibinfo
  {year} {1998})\BibitemShut {NoStop}%
  %
  %
  %
\bibitem [{\citenamefont {Palatini}(1933)}]{zbMATH03012628}%
  \BibitemOpen
  \bibfield  {author} {\bibinfo {author} {\bibfnamefont {A.}~\bibnamefont
  {Palatini}},\ }
%  
   % 
  {\enquote {\bibinfo {title} {{Concetto di vettore generalizzato prodotto interno, prodotto esterno, divergenza e rotore. Teoremi generali della divergenza, del rotore e di Stokes}},}\ }
%  
%  
  \href {http://www.numdam.org/item?id=RSMUP_1933__4__122_0}
  {\bibfield  {journal} {\bibinfo  {journal} {Rend. Semin. Mat. Univ. Padova}\
  }\textbf {\bibinfo {volume} {4}},\ \bibinfo {pages} {122-139} (\bibinfo {year}
  {1933})}\BibitemShut {NoStop}%
  %
  %
\bibitem [{\citenamefont
  {Silagadze}(2002{\natexlab{a}})}]{Silagadze:Maxwell_7D}%
  \BibitemOpen
  \bibfield  {author} {\bibinfo {author} {\bibfnamefont {Z.~K.}\ \bibnamefont
  {Silagadze}},\ }
  %  
   % 
  {\enquote {\bibinfo {title} {{Feynman's derivation of Maxwell equations and extra dimensions}},}\ }
  %
  %
  \href {http://arxiv.org/abs/hep-ph/0106235} {\bibfield
  {journal} {\bibinfo  {journal} {Ann. Fond. Louis de Broglie}\ }\textbf
  {\bibinfo {volume} {27}},\ \bibinfo {pages} {241-255} (\bibinfo {year}
  {2002}{\natexlab{a}})},\ \bibinfo {note} {{Special issue on contemporary
  electrodynamics}}\BibitemShut {NoStop}%
  %
  %
  %
\bibitem [{\citenamefont
  {Silagadze}(2002{\natexlab{b}})}]{Silagadze:multi_dim_vector_prod}%
  \BibitemOpen
  \bibfield  {author} {\bibinfo {author} {\bibfnamefont {Z.~K.}\ \bibnamefont
  {Silagadze}},\ }
  %
     % 
  {\enquote {\bibinfo {title} {{Multi-dimensional vector product}},}\ }
  %
  %
  \href {http://arxiv.org/abs/math/0204357} {\bibfield
  {journal} {\bibinfo  {journal} {J. Phys. A: Math. Gen.}\ }\textbf {\bibinfo
  {volume} {35}},\ \bibinfo {pages} {4949-4953} (\bibinfo {year}
  {2002}{\natexlab{b}})}\BibitemShut {NoStop}%
%-----------------------------------------------------------
%
\end{thebibliography}
%Termine documento.
\end{document}